\documentclass[12pt,oneside]{amsart}
\usepackage{geometry}                
\usepackage{graphicx}
\usepackage{amssymb}
\usepackage{epstopdf}
\usepackage{amsmath,amsthm,amscd,amssymb}
\usepackage{latexsym}
\usepackage[colorlinks,citecolor=red,pagebackref,hypertexnames=false]{hyperref}

\numberwithin{equation}{section}

\theoremstyle{plain}
\newtheorem{theorem}{Theorem}[section]
\newtheorem{lemma}[theorem]{Lemma}
\newtheorem{corollary}[theorem]{Corollary}

\theoremstyle{definition}
\newtheorem{definition}[theorem]{Definition}

\theoremstyle{remark}
\newtheorem{remark}[theorem]{Remark}

\newtheorem{case[theorem]}{Case}

\date{August 22, 2011}      
\author{Alex Iosevich, Mihalis Mourgoglou and Eyvindur Ari Palsson}
\address{Department of Mathematics, University of Rochester, Rochester, NY}
\email{iosevich@math.rochester.edu}
\address{D\'{e}partement de Math\'{e}matiques\\ UMR 8628 Universit\'{e} Paris-Sud 11-CNRS \\  B\^atiment 425\\ Facult\'{e} des Sciences d'Orsay\\ Universit\'{e} Paris-Sud 11\\ F-91405 Orsay Cedex}
\email{mihalis.mourgoglou@math.u-psud.fr}
\address{Department of Mathematics, University of Rochester, Rochester, NY}
\email{palsson@math.rochester.edu}

\thanks{This work was partially supported by the NSF Grant DMS10-45404.}

\title{\parbox{14cm}{\centering{On angles determined by fractal subsets of the Euclidean space via Sobolev bounds for bi-linear operators}}}

\begin{document}


\begin{abstract} We prove that if the Hausdorff dimension of a compact subset of ${\mathbb R}^d$ is greater than $\frac{d+1}{2}$, then the set of angles determined by triples of points from this set has positive Lebesgue measure. Sobolev bounds for bi-linear analogs of generalized Radon transforms and the method of stationary phase play a key role. These results complement those of V. Harangi, T. Keleti, G. Kiss, P. Maga, P. Mattila and B. Stenner in (\cite{HKKMMS10}). We also obtain new upper bounds for the number of times an angle can occur among $N$ points in ${\mathbb R}^d$, $d \ge 4$, motivated by the results of Apfelbaum and Sharir (\cite{AS05}) and Pach and Sharir (\cite{PS92}).  We then use this result to establish sharpness results in the continuous setting. Another sharpness result relies on the distribution of lattice points on large spheres in higher dimensions. 
\end{abstract}  

\maketitle


\section{Introduction} 

In this paper we study angles determined by subsets of the Euclidean space of a given Hausdorff dimension. 

\begin{definition} Given $E \subset {\mathbb R}^d$, $d \ge 2$, let $\theta(x^1,x^2,x^3)$ denote the interior angle of the triangle with vertices at $x^1x^2x^3$ , at $x^1$, where $x^j \in E$. Define the angle set 
$${\mathcal A}(E)=\{\theta(x^1,x^2,x^3): x^j \in E\}.$$ \end{definition} 

The question we ask is, how large does the Hausdorff dimension of $E \subset {\mathbb R}^d$, $d \ge 2$, need to be to ensure that the Lebesgue measure of ${\mathcal A}(E)$ is positive. Similarly, we would like to know if the angles are uniformly distributed in the sense that a small neighborhood of a given angle does not arise more often than is its share. 

Our results are partly motivated by a paper due to V. Harangi, T. Keleti, G. Kiss, P. Maga, P. Mattila and B. Stenner (\cite{HKKMMS10}) where it is proved that if the Hausdorff dimension of $E \subset {\mathbb R}^d$, $d \ge 2$ is greater than $d-1$, then every angle $\theta \in [0, \pi]$ is in ${\mathcal A}(E)$. The authors also prove that if the Hausdorff dimension is greater than $\frac{d}{2}$, if $d$ is even, and $\frac{d+1}{2}$ if $d$ is odd, then the angle $\frac{\pi}{2}$ is in ${\mathcal A}(E)$. Furthermore, they demonstrate that the threshold $d-1$ is best possible for $\theta=\pi$. 

In this paper we show that the Hausdorff dimensional threshold $\frac{d+1}{2}$ is sufficient to ensure that ${\mathcal A}(E)$ has positive Lebesgue measure. This will follow from the fact that in this regime, no angle is overrepresented. 

\begin{definition} Let $E \subset {\mathbb R}^d$, $d \ge 2$. We say that an angle $\alpha \in [0, \pi]$ is {\it equitably represented} in $\mathcal{A}(E)$ if for every Frostman measure $\mu$ supported on $E$ and any $\epsilon>0$, 
\begin{equation}\label{def}
\mu \times \mu \times \mu \{(x^1,x^2,x^3): \alpha-\epsilon \leq \theta(x^1,x^2,x^3) \leq \alpha+\epsilon \} \lesssim  \epsilon .
\end{equation}
Here and throughout, $X \lesssim Y$ means that there exists $C>0$ such that $X \leq CY$. 
\end{definition}

Recall that a probability measure $\mu$ on a compact set $E\subset\mathbb{R}^d$ is a {\it Frostman measure}  if, for any ball $B_{\delta}$ of radius $\delta$, 
\begin{equation} \label{frostman} \mu(B_{\delta}) \lessapprox \delta^s, \end{equation} 
where $s=dim_{{\mathcal H}}(E))$. For discussion and  proof of the existence of such measures see, e.g., \cite{M95}. 

Our main result is the following. 

\begin{theorem} \label{main} Let $E$ be a compact subset of ${\mathbb R}^d$ of Hausdorff dimension greater than $\frac{d+1}{2}$. Then every $\alpha \in [0, \pi]$ is equitably represented in $\mathcal{A}(E)$. 
\end{theorem}

\begin{remark} We note that the implicit constant in (\ref{def}) that we obtain depends only on the dimension $d$ and on $I_s(\mu)$ for some $s \in (\frac{d+1}{2}, dim_{{\mathcal H}}(E))$, where
$$ I_s(\mu)=\int \int {|x-y|}^{-s} d\mu(x) d\mu(y).$$
\end{remark}

\begin{corollary} \label{maincorollary}  Let $E$ be a compact subset of ${\mathbb R}^d$ of Hausdorff dimension greater than $\frac{d+1}{2}$. Then the Lebesgue measure of ${\mathcal A}(E)$ is positive. 
\end{corollary} 

The key tool is the following bi-linear estimate. 
\begin{theorem} \label{bilinearangles} Let 
$$T_{\mu_t}(f,g)(x)=\int \int f(x-u) g(x-v) d\mu_t(u,v),$$
and
$$\tilde{T}_{\mu_t}(f,g)(x)=\int \int f(x+u-v) g(x+u) d\mu_t(u,v),$$
where $\mu_t$ is the Leray measure on the set 
$$ \left\{(u,v) \in B \times B: \frac{u}{|u|} \cdot \frac{v}{|v|}=t \right\},$$ where $-1 \leq t \leq 1$ and $B$ denotes the unit ball. Then 
\begin{equation} \label{sobolev} \|T_{\mu_t}(f,g)\|_{L^1({\mathbb R}^d)} \lesssim \| f \|_{ L^2_{-\alpha}({\mathbb R}^d)} \| g \|_{L^2_{-\beta}({\mathbb R}^d)} \end{equation}
and
\begin{equation} \label{sobolev2} \|\tilde{T}_{\mu_t}(f,g)\|_{L^1({\mathbb R}^d)} \lesssim \| f \|_{ L^2_{-\alpha}({\mathbb R}^d)} \| g \|_{L^2_{-\beta}({\mathbb R}^d)} \end{equation}
with $\alpha, \beta \ge 0$, $\alpha+\beta=d-1$, for all non-negative functions $f$ and $g$. 
\end{theorem} 

\begin{remark} Here and throughout $L^2_{-\alpha}({\mathbb R}^d)$ denotes the Sobolev space with negative index where 
$$ {||f||}^2_{L^2_{-\alpha}({\mathbb R}^d)}=\int_{{\mathbb R}^d} {|\widehat{f}(\xi)|}^2 {|\xi|}^{-2\alpha} d\xi$$ and 
$$ \widehat{f}(\xi)=\int_{{\mathbb R}^d} e^{-2 \pi i x \cdot \xi} f(x) dx,$$ the Fourier transform of $f: {\Bbb R}^d \to {\Bbb R}$. 
\end{remark} 

Theorem \ref{bilinearangles} will be deduced from the following two results inspired by a more general result proved by Allan Greenleaf and the first and third listed authors in \cite{GIP11}.
\begin{theorem} \label{miracle} Let 
$$ T_{\mu}(f,g)(x)=\int \int f(x-u) g(x-v) d\mu(u,v),$$ where $d\mu$ is a positive finite measure on ${\mathbb R}^d$. Suppose there exists $C>0$ such that 
\begin{equation} \label{decay} |\widehat{\mu}(\xi, \eta)| \leq C({|\xi|+|\eta|)}^{-\gamma} \end{equation} for some $\gamma>0$. Then 
\begin{equation} \label{orgasm}  \|T_{\mu_t}(f,g)\|_{L^1({\mathbb R}^d)} \lesssim \| f \|_{ L^2_{-\alpha}({\mathbb R}^d)} \| g \|_{L^2_{-\beta}({\mathbb R}^d)} \end{equation} with $\alpha, \beta \ge 0$, $\alpha+\beta=\gamma$, for all non-negative functions $f$ and $g$. \end{theorem} 

\begin{theorem} \label{miracle2} Let 
$$ \tilde{T}_{\mu}(f,g)(x)=\int \int f(x+u-v) g(x+u) d\mu(u,v),$$ where $d\mu$ is a positive finite measure on ${\mathbb R}^d$. Suppose there exists $C>0$ such that 
\begin{equation} \label{decay2} |\widehat{\mu}(\xi, \eta)| \leq C({|\xi|+|\eta|)}^{-\gamma} \end{equation} for some $\gamma>0$. Then 
\begin{equation} \label{orgasm2}  \|T_{\mu_t}(f,g)\|_{L^1({\mathbb R}^d)} \lesssim \| f \|_{ L^2_{-\alpha}({\mathbb R}^d)} \| g \|_{L^2_{-\beta}({\mathbb R}^d)} \end{equation} with $\alpha, \beta \ge 0$, $\alpha+\beta=\gamma$, for all non-negative functions $f$ and $g$. \end{theorem} 

\begin{remark} Going carefully through the proof of \ref{miracle} one can see that the condition \ref{decay} for all $\xi$, $\eta$ can be significantly weakened. Roughly speaking, we only need the decay rate to hold near the hyperplane $\xi+\eta=(0, \dots, 0)$. This came in quite handy in \cite{GI11}. A more general theorem can be found in \cite{GIP11}. Similar comments also apply to the proof of \ref{miracle2}.
\end{remark} 

\vskip.125in 

\subsection{Sharpness of results} 

We show that a necessary lower bound on the Hausdorff dimension of $E \subset {\mathbb R}^d$, $d \ge 2$, is $\frac{d}{2}$ in order to ensure that every angle is not overrepresented. Thus it is an open question whether our bound $\frac{d+1}{2}$ from the main theorem can be improved or whether an example can be found that shows the sharpness of our result.

\begin{theorem}\label{notsharp} For every $d \ge 2$ and $s \in (0, \frac{d}{2})$, there exists $E \subset {\mathbb R}^d$ of Hausdorff dimension $s$ such that $\frac{\pi}{2}$ is not equitably represented in $\mathcal{A}(E)$. \end{theorem}

The main ingredient in the proof is the following generalization to $\mathbb{R}^d$ of a theorem by Apfelbaum and Sharir in \cite{AS05}, which they state in $\mathbb{R}^3$.

\begin{theorem}\label{counting}
Let $P_n=\lbrace 1,\ldots , \lfloor n^{1/d} \rfloor \rbrace^d$. Then the number of triplets $(p,q,r) \in P_n^3$ such that $\angle pqr = \pi/2$ is $\Omega(n^{3-\frac{2}{d}})$. Here and throughout, $X=\Omega(Y)$ with the controlling parameter $n$ means that there exists $c>0$ such that $X \ge cY$ with $c$ independent of $n$. 
\end{theorem}

\begin{remark} \label{donotknowshit} We do not know if a version of Theorem \ref{notsharp} exists for angles other than $\frac{\pi}{2}$ in the same range of exponents. Similarly, it is not known if a version of Theorem \ref{counting} exists for other angles. See \cite{AS05} for a detailed discussion of this issue. 
\end{remark} 

We also have a necessity result with respect to the positive Lebesgue measure of ${\mathcal A}(E)$, the set of angles. Since $E$ may be a subset of a line, we immediately see that in order to ensure that the Lebesgue measure of ${\mathcal A}(E)$ is positive, we must assume that the Hausdorff dimension of $E$ is greater than one. Another example is the following. 

\begin{theorem} \label{notsharp2} For every $d \ge 4$ and every $s<\frac{d-2}{2}$, there exists a sequence of sets $E_n$ with the following properties: 

\begin{itemize} 
\item Each $E_n$ is a union of balls. 

\item Each $E_n$ uniformly $s$-dimensional in the sense that if $\mu_n$ is the normalized Lebesgue measure on $E_n$, then 
\begin{equation} \label{energy} \int \int {|x-y|}^{-s} d\mu_n(x) d\mu_n(y) \approx 1. \end{equation} 

\item The Lebesgue measure of ${\mathcal A}(E_n)$ tends to $0$ as $n$ tends to infinity. 

\end{itemize} 
\end{theorem} 

\begin{remark} It would be interesting to sharpen up Theorem \ref{notsharp2} to say that for every $s<\frac{d-2}{2}$, $d \ge 4$, there exists $E \subset {\Bbb R}^d$ of Hausdorff dimension $s$ such that the Lebesgue measure of ${\mathcal A}(E)$ is $0$. 
\end{remark}

\vskip.125in 

\subsection{Applications to discrete geometry} 
\label{discreteapplications}

The following results were obtained by Pach and Sharir, in \cite{PS92}, and Apfelbaum and Sharir, in \cite{AS05}. In \cite{PS92}, it is shown that for a set of $n$ points in $\mathbb{R}^2$, no angle can occur more than $cn^2 \log n$ times. Since there are about $n^3$ triples of points, this implies that there must be at least $c\frac{n}{\log n}$ distinct angles. In \cite{AS05}, it is shown that for a set of $n$ points in $\mathbb{R}^3$, no angle can occur more than $cn^{\frac{7}{3}}$ times, which gives a lower bound of at least $cn^{\frac{2}{3}}$ distinct angles. They also show that for a set of $n$ points in $\mathbb{R}^4$, no angle besides $\frac{\pi}{2}$ can occur more than $cn^{\frac{5}{2}}\beta(n)$ times, where $\beta(n)$ grows extremely slowly with respect to $n$. This means that there must be about $n^{\frac{1}{2}}(\beta(n))^{-1}$ distinct angles.

In dimensions four and higher, no results are currently available. In order to describe our main result in this direction, we need the following definition. 
\begin{definition} \label{adaptablemama} Let $P$ be a set of $n$ points contained in ${[0,1]}^d$, $d \ge 2$. Define the measure
$$ d \mu^s_P(x)=n^{-1} \cdot n^{\frac{d}{s}} \cdot \sum_{p \in P} \chi_{B_{n^{-\frac{1}{s}}}(p)}(x)dx,$$ where $\chi_{B_{n^{-\frac{1}{s}}}(p)}(x)$ is the characteristic function of the ball of radius $n^{-\frac{1}{s}}$ centered at $p$. 

We say that $P$ is \underline{$s$-adaptable} if $P$ is $n^{-\frac{1}{s}}$-separated and 

$$ I_s(\mu_P)=\int \int {|x-y|}^{-s} d\mu^s_P(x) d\mu^s_P(y)<\infty.$$ 

\vskip.125in 

This is equivalent to the statement 

$$ n^{-2} \sum_{p \not=p' \in P} {|p-p'|}^{-s} \lesssim 1.$$ 

\end{definition}

\vskip.125in 

To put it simply, $s$-adaptability means that a discrete point set $P$ can be thickened into a set which is uniformly $s$-dimensional in the sense that its energy integral of order $s$ is finite. Unfortunately, it is shown in \cite{IRU10} that there exist finite point sets which are not $s$-adaptable for certain ranges of the parameter $s$. However, many commonly used classes of discrete sets, such as homogeneous sets, studied, for example, by Solymosi and Vu, are indeed $s$-adaptable for $0 \leq s \leq d$. See \cite{IRU10} for a detailed description of these issues. Our main discrete geometric result is the following. 
\begin{theorem} \label{discrete}
Let $P \subset {\mathbb R}^d$, $\#P=N$, $d \ge 2$, be an $s$-adaptable set for $s>\frac{d+1}{2}$. Then 
$$ \# \{(x^1,x^2,x^3) \in P \times P \times P: \theta(x^1,x^2,x^3)=\theta_0 \} \lesssim N^{3-\frac{1}{s}}.$$ 
\end{theorem} 

In dimensions two and three, these exponents are not as good as the aforementioned results of Apfelbaum and Sharir and Pach and Sharir. However, Theorem \ref{discrete} gives non-trivial exponents in all dimensions. 

\vskip.125in 

\section{Proof of Corollary \ref{maincorollary}} 

\vskip.125in 

Define 
$$ \nu^{\epsilon}(t)=\frac{1}{\epsilon} \mu \times \mu \times \mu \left\{(x,y,z): t - \epsilon \leq \frac{x-z}{|x-z|} \cdot \frac{y-z}{|y-z|} \leq t+\epsilon \right\},$$ where $\mu$ is a Borel measure on $E$. Theorem  \ref{main} states precisely that 
\begin{equation} \label{erection} \nu^{\epsilon}(t) \lesssim 1. \end{equation} 

Cover ${\mathcal A}(E)$ by $\cup_i (t_i-\epsilon_i, t_i+\epsilon_i)$. It follows that 
$$ 1=\mu \times \mu \times \mu \{E \times E \times E \}$$
$$ \leq \sum_i \mu \times \mu \times \mu \left\{(x,y,z):  t_i-\epsilon_i \leq \frac{x-z}{|x-z|} \cdot \frac{y-z}{|y-z|} \leq t_i+\epsilon_i \right\} \lesssim \sum_i \epsilon_i,$$ where the last inequality follows by (\ref{erection}). We conclude that $\sum_i \epsilon_i \gtrsim 1$, which implies that the Lebesgue measure of ${\mathcal A}(E)$ is positive.

\vskip.125in 

\section{Proof of Theorem \ref{main}} 

\vskip.125in 

Let $\mu$ be a Frostman measure, supported on $E \subset {\mathbb R}^d$. Let $\rho$ be a non-negative, smooth function, equal to 1 on $[-\frac{1}{4}, \frac{1}{4}]$, supported in $[-1,1]$ with $\| \rho \|_{L^1(\mathbb{R}^d)} = 1$. Set $\mu^{\delta} = \mu * \rho_{\delta}$ where $\rho_{\delta}(x) = \delta^{-d}\rho(\frac{x}{\delta})$, which means that $\mu^{\delta}$ is a smooth approximation of $\mu$ and tends to $\mu$ as $\delta$ tends to $0$. We establish a bound
$$ \frac{1}{\epsilon}\ \mu^{\delta} \times \mu^{\delta} \times \mu^{\delta} \left\{(x,y,z): t-\epsilon \leq \frac{x-z}{|x-z|} \cdot 
\frac{y-z}{|y-z|} \leq t+\epsilon \right\} \lesssim 1$$
independent of $\delta$ and thus by passing to the limit we establish our theorem.

Write
$$ \frac{1}{\epsilon}\ \mu^{\delta} \times \mu^{\delta} \times \mu^{\delta} \left\{(x,y,z): t-\epsilon \leq \frac{x-z}{|x-z|} \cdot 
\frac{y-z}{|y-z|} \leq t+\epsilon \right\}$$ 
$$=\epsilon^{-1}\ \int \int \int_{\left\{(x,y): t-\epsilon \leq \frac{x-z}{|x-z|} \cdot \frac{y-z}{|y-z|} \leq t+\epsilon \right\}} d\mu^{\delta}(x) d\mu^{\delta}(y) d\mu^{\delta}(z)$$
\begin{equation} \label{foreplay}=\langle T^{\epsilon}(\mu^{\delta}, \mu^{\delta}), \mu^{\delta} \rangle ,\end{equation} where 
$$ T^{\epsilon}(f,g)(z)=\epsilon^{-1} \int \int_{\left\{(x,y): t-\epsilon \leq \frac{x-z}{|x-z|} \cdot \frac{y-z}{|y-z|} \leq t+\epsilon \right\}} f(x)g(y) dxdy$$
$$ \equiv \int \int f(z-u) g(z-v) d\mu_t^{\epsilon}(u,v)$$
and $\langle \cdot,\cdot \rangle$ is the $L^2(\mathbb{R}^d)$ inner product.

Now define
$$ F(\alpha) := \langle T^{\epsilon}(\mu^{\delta}_{-\alpha}, \mu^{\delta}), \mu^{\delta}_{\alpha} \rangle $$
where
\begin{equation}\label{spongebob}
\mu^{\delta}_{\beta}(x) := \frac{2^{\frac{d-\beta}{2}}}{\Gamma\left(\beta/2 \right)}(\mu^{\delta} * |\cdot |^{-d+\beta})(x)
\end{equation}
initially defined for $\text{Re}(\beta) > 0$, is extended to the complex plane by analytic continuation. Since $\mu^{\delta}_{\beta}$ is smooth and we are in a compact setting then we have trivial bounds on $F(\alpha) = \langle T^{\epsilon}(\mu^{\delta}_{-\alpha}, \mu^{\delta}), \mu^{\delta}_{\alpha} \rangle$ with constants depending on $\delta$. Observe that $\widehat{\mu^{\delta}_{\beta}} (\xi) = C_{\beta,d}\widehat{\mu}(\xi)\widehat{\rho}(\delta\xi)|\xi|^{-\beta}$ where
\begin{equation}\label{blowup}
C_{\beta,d} = \frac{(2\pi)^{\frac{d}{2}}2^{\frac{\beta}{2}}}{\Gamma(\frac{d-\beta}{2})}.
\end{equation}
See, e.g., page 192 in \cite{GS58} for related calculations. By Plancherel then $\mu_{\beta}^\delta$ is an $L^2({\mathbb R}^d)$ function with bounds depending on $\delta$. Taking the modulus in (\ref{spongebob}), we see that
$$|\mu_\beta^\delta(x)|\leq \left| \frac{2^{\frac{d-\beta}{2}}}{\Gamma\left(\beta/2 \right)} \right| (\mu^{\delta} * |\cdot |^{-d+\text{Re}(\beta)})(x) = 2^{\frac{d-\text{Re}(\beta)}{2}}\frac{\Gamma(\text{Re}(\beta)/2)}{|\Gamma(\beta/2)|}\mu_{\text{Re}(\beta)}^\delta(x)$$
and note that the right hand side is non-negative.

We now estimate
$$ |F(\alpha)| \leq  {||T^{\epsilon}(\mu^{\delta}_{-\alpha}, \mu^{\delta})||}_1 \cdot {||\mu^{\delta}_{\alpha}||}_{\infty} \leq {\left|\left|T^{\epsilon}\left(2^{\frac{d-\text{Re}(\alpha)}{2}}\frac{\Gamma(\text{Re}(\alpha)/2)}{|\Gamma(\alpha/2)|}\mu^{\delta}_{-\text{Re}(\alpha)}, \mu^{\delta}\right)\right|\right|}_1 \cdot {||\mu^{\delta}_{\text{Re}(\alpha)}||}_{\infty} .$$
Using \ref{sobolev} from theorem \ref{bilinearangles} we have
\begin{multline*}
{\left|\left|T^{\epsilon}\left(2^{\frac{d-\text{Re}(\alpha)}{2}}\frac{\Gamma(\text{Re}(\alpha)/2)}{|\Gamma(\alpha/2)|}\mu^{\delta}_{-\text{Re}(\alpha)}, \mu^{\delta}\right)\right|\right|}_1 \\
 \lesssim {\left|\left|2^{\frac{d-\text{Re}(\alpha)}{2}}\frac{\Gamma(\text{Re}(\alpha)/2)}{|\Gamma(\alpha/2)|}\mu^{\delta}_{-\text{Re}(\alpha)}\right|\right|}_{ L^2_{-\frac{3(d-1)}{4}}({\mathbb R}^d)}  \cdot {||\mu^{\delta}||}_{ L^2_{-\frac{(d-1)}{4}}({\mathbb R}^d)}
\end{multline*}
so we can conclude
$$ |F(\alpha)| \lesssim {||\mu^{\delta}_{-\text{Re}(\alpha)}||}_{ L^2_{-\frac{3(d-1)}{4}}({\mathbb R}^d)}  \cdot {||\mu^{\delta}||}_{ L^2_{-\frac{(d-1)}{4}}({\mathbb R}^d)}\cdot{||\mu^{\delta}_{\text{Re}(\alpha)}||}_{\infty} $$
where our implicit constants depend on gamma functions.
A standard calculation shows
$$ {||\mu^{\delta}||}_{ L^2_{-\frac{(d-1)}{4}}({\mathbb R}^d)} \lesssim \left(I_{d-2\cdot \frac{d-1}{4}}(\mu^{\delta})\right)^{1/2} = \left(I_{\frac{d+1}{2}}(\mu^{\delta})\right)^{1/2}$$
where
$$ I_s(\mu)=\int \int {|x-y|}^{-s} d\mu(x) d\mu(y).$$

Now take $\alpha$ with $\text{Re}(\alpha) = \frac{d-1}{2}$. We then have
\begin{align*}
|F(\alpha)| &\lesssim {||\mu^{\delta}_{-\frac{d-1}{2}}||}_{ L^2_{-\frac{3(d-1)}{4}}({\mathbb R}^d)}  \cdot {||\mu^{\delta}||}_{ L^2_{-\frac{(d-1)}{4}}({\mathbb R}^d)}\cdot{||\mu^{\delta}_{\frac{d-1}{2}}||}_{\infty}\\
&\lesssim {||\mu^{\delta}||}_{ L^2_{-\frac{(d-1)}{4}}({\mathbb R}^d)}\cdot {||\mu^{\delta}||}_{ L^2_{-\frac{(d-1)}{4}}({\mathbb R}^d)} \cdot {||\mu^{\delta}_{\frac{d-1}{2}}||}_{\infty}\\
&\lesssim I_{\frac{d+1}{2}}(\mu^{\delta}) \cdot {||\mu^{\delta}_{\frac{d-1}{2}}||}_{\infty}.
\end{align*}
and since $\mu$ is supported on a set of Hausdorff dimension greater than $\frac{d+1}{2}$ we can bound
$$ I_{\frac{d+1}{2}}(\mu^{\delta}) \lesssim 1 $$
with a bound independent of $\delta$. Note that in the above calculations we can not allow $-\frac{d-1}{4}$ to be a negative integer because the implicit constants in our bounds depend on the gamma function, $\Gamma$, evaluated at this point. We can however always consider instead $\alpha$ with $\text{Re}(\alpha) = \frac{d-1}{2} - \kappa$ where $0<\kappa \ll dim_{{\mathcal H}}(E) - \frac{d+1}{2}$ since $\frac{d+1}{2}<dim_{{\mathcal H}}(E)$ and we would obtain similar bounds and note that we still have $I_{\frac{d+1}{2}+\kappa}(\mu^{\delta}) \lesssim 1$.

To show that $|F(\alpha)| \lesssim 1$ for $\alpha$ with $\text{Re}(\alpha) = \frac{d-1}{2}$ we are thus left with showing
$$ {||\mu^{\delta}_{\frac{d-1}{2}}||}_{\infty} \lesssim 1 .$$
This however follows from the following calculations
$$ |\mu^{\delta}_{\frac{d-1}{2}}| \lesssim \int |x-y|^{-d+\frac{d-1}{2}}d\mu^{\delta}(y) \approx \sum_m 2^{m\frac{d+1}{2}}\int_{|x-y|\approx 2^{-m}}d\mu^{\delta}(y) \lesssim \sum_m 2^{m\frac{d+1}{2}} 2^{-ms} \lesssim 1  ,$$
since $\mu$ is a Frostman measure on a set of Hausdorff dimension $s > \frac{d+1}{2}$.

Note that $F(\alpha)$ is analytic and bounded on $[-\frac{d-1}{2}, \frac{d-1}{2}]$. We have shown that $F(\alpha)\lesssim 1$ if $\text{Re}(\alpha) = \frac{d-1}{2}$ with a bound independent of $\delta$. If we could show a similar bound for $\text{Re}(\alpha) =- \frac{d-1}{2}$ then by a generalization of the three lines lemma, due to I. I. Hirschman \cite{H52}, we would have $F(z)\lesssim 1$, with a bound independent of $\delta$, for all $z$ with $- \frac{d-1}{2} \leq \text{Re}(z) \leq \frac{d-1}{2}$. In particular we would have $F(0) \lesssim 1$ which would prove our theorem.

In \cite{GI11} Iosevich and Greenleaf used a similar strategy where they obtained $F(\alpha)\lesssim 1$ and then immediately by symmetry $F(-\alpha)\lesssim 1$. Here we do not immediately obtain such a bound by symmetry. We must go back to our trilinear form, rewrite it and use results on the boundedness of an operator which is an adjoint to the bilinear operator we considered above.

By changing variables a couple of times we obtain
\begin{align*}
F(-\alpha) &=\epsilon^{-1}\ \int \int \int_{\left\{(x,y): t-\epsilon \leq \frac{x-z}{|x-z|} \cdot \frac{y-z}{|y-z|} \leq t+\epsilon \right\}} \mu_{\alpha}^{\delta}(x) \mu^{\delta}(y) \mu_{-\alpha}^{\delta}(z) dxdydz\\
&= \epsilon^{-1}\ \int \int \int_{\left\{(u,v): t-\epsilon \leq \frac{u}{|u|} \cdot \frac{v}{|v|} \leq t+\epsilon \right\}} \mu_{\alpha}^{\delta}(z-u) \mu^{\delta}(z-v) \mu_{-\alpha}^{\delta}(z) dudvdz\\
&= \epsilon^{-1}\ \int \int \int_{\left\{(u,v): t-\epsilon \leq \frac{u}{|u|} \cdot \frac{v}{|v|} \leq t+\epsilon \right\}} \mu_{\alpha}^{\delta}(z) \mu^{\delta}(z+u-v) \mu_{-\alpha}^{\delta}(z+u) dudvdz\\
&=\langle \tilde{T}^{\epsilon}(\mu_{-\alpha}^{\delta}, \mu^{\delta}), \mu_{\alpha}^{\delta} \rangle
\end{align*}
where 
$$ \tilde{T}^{\epsilon}(f,g)(z)=\epsilon^{-1} \int \int_{\left\{(u,v): t-\epsilon \leq \frac{u}{|u|} \cdot \frac{v}{|v|} \leq t+\epsilon \right\}} f(z+u-v)g(z+u) dudv$$
$$ \equiv \int \int f(z+u-v) g(z+u) d\mu_t^{\epsilon}(u,v)$$

After this rewrite it is clear that the same procedure as above goes through unchanged, except we must use \ref{sobolev2} from theorem \ref{bilinearangles} to estimate the bilinear operator.

Thus the proof of Theorem \ref{main} has been reduced to proving Theorem \ref{bilinearangles} and this is where we now turn our attention.

\vskip.25in 

\section{Proof of Theorem \ref{bilinearangles}} 

As we point out in the introduction, Theorem \ref{bilinearangles} would follow from Theorem \ref{miracle} and Theorem \ref{miracle2} if we could show that (\ref{decay}) holds with $\gamma=d-1$. Suppose that $|\eta| \ge c|\xi|$. Then parameterize the set 
$$ \left\{(u,v) \in {[0,1]}^d \times {[0,1]}^d: \frac{u}{|u|} \cdot \frac{v}{|v|}=t \right\}$$ as 
\begin{equation} \label{bendover} \left\{(u, \lambda \theta u): \lambda \in [0, \sqrt{d}]: u \in {[0,1]}^d, \theta \in \Omega_{t, \frac{u}{|u|}} \right\},\end{equation} where 
$$ \Omega_{t, \frac{u}{|u|}}=\{ \theta \in SO(d): \theta u \cdot u={|u|}^2 t \}.$$ 

It follows that 
\begin{equation} \label{bendover2} \widehat{\mu}_t(\xi, \eta)=\iiint e^{-2 \pi i (u \cdot \xi+ \lambda \theta u \cdot \eta)} \psi(|u|) \psi_0(\lambda) d\Omega_{t, \frac{u}{|u|}}(\theta)du d\lambda, \end{equation} where $d\Omega_{t, \frac{u}{|u|}}(\theta)$ is the restriction of the Haar measure on $SO(d)$ to $\Omega_{t, \frac{u}{|u|}}$ and $\psi, \psi_0$ are smooth cut-off functions. The modulus of the expression (\ref{bendover2}) equals 
$$ \left|\iint e^{-2 \pi i u \cdot \xi} \widehat{\psi}_0(\theta u \cdot \eta) \psi(|u|) d\Omega_{t, \frac{u}{|u|}}(\theta) du \right| $$
$$ \leq \iint |\widehat{\psi}_0(\theta u \cdot \eta)| \psi(|u|) d\Omega_{t, \frac{u}{|u|}}(\theta)du$$
$$ \lesssim \sup\limits_{\omega\in S^{d-1}} \iint  |\widehat{\psi}_0(\theta u \cdot \eta)| \psi(|u|) d\Omega_{t,\omega}(\theta)du$$
$$ \lesssim  \sup\limits_{\omega\in S^{d-1}} \iint \sum\limits_{k=1}^{\infty}2^{-k+1}|\lbrace k-1\lesssim |\theta u \cdot \eta| \lesssim k\rbrace|\psi(|u|) d\Omega_{t,\omega}(\theta)du$$
$$ = \sup\limits_{\omega\in S^{d-1}} \iint \sum\limits_{k=1}^{\infty}2^{-k+1}\left|\left\lbrace\frac{k-1}{|\eta|} \lesssim \left| u \cdot \frac{\theta^{T}\eta}{|\eta|}\right| \lesssim \frac{k}{|\eta|}\right\rbrace\right|\psi(|u|) du d\Omega_{t,\omega}(\theta)$$
$$ \lesssim \sup\limits_{\omega\in S^{d-1}} \iint \sum\limits_{k=1}^{\infty}2^{-k+1}\left|\left\lbrace\frac{k-1}{|\eta|} \lesssim \left| r\sigma \cdot \frac{\theta^{T}\eta}{|\eta|}\right| \lesssim \frac{k}{|\eta|}\right\rbrace\right|\psi(r) r^{d-1} dr d\sigma d\Omega_{t,\omega}(\theta)$$
$$ \lesssim \sup\limits_{\omega\in S^{d-1}} \iint \sum\limits_{k=1}^{\infty}2^{-k+1}k\frac{1}{r^{d-1}|\eta|^{d-1}}\psi(r) r^{d-1}dr d\Omega_{t,\omega}(\theta)$$
$$ \lesssim \frac{1}{|\eta|^{d-1}} $$
$$ \lesssim \frac{1}{(|\xi|+|\eta|)^{d-1}} $$

If $|\eta| < c|\xi|$ then run the same argument as above, except with the roles of $\xi$ and $\eta$ reversed, and obtain the same bound.

\vskip.125in 

\section{Proof of Theorem \ref{miracle}} 
\label{miraclesection} 

In this section we prove Theorem \ref{miracle}. In order to study the integrability of
$$ T_{\mu}(f,g)(x)=\iint f(x-u) g(x-v) d\mu(u,v)$$
we look at
$$ \int \left|T_{\mu}(f,g)(x)\psi\left(\frac{x}{R}\right)\right| dx $$
where $\psi$ is a non-negative bump function with support in $[-2,2]$ and identically equal to $1$ on $[-1,1]$. We obtain bounds independent on $R$ and thus Theorem \ref{miracle} follows by a standard limiting argument.

Assum $f$ and $g$ are positive functions so that we can write
$$ \int T_{\mu}(f,g)(x)\psi\left(\frac{x}{R}\right) dx .$$
Viewing the above expression from the Fourier side we obtain
$$ R^d \iint \widehat{f}(\xi)\widehat{g}(\eta) \widehat{\mu}(\xi,\eta)\widehat{\psi}(R(\xi + \eta)) d\xi d\eta .$$
By using estimate \ref{decay} we can bound the above expression by
$$ R^d \iint |\widehat{f}(\xi)| |\widehat{g}(\eta)|C(|\xi|+|\eta|)^{-\gamma}|\widehat{\psi}(R(\xi + \eta))| d\xi d\eta $$
$$ \lesssim R^d \iint |\widehat{f}(\xi)| |\widehat{g}(\eta)| |\xi|^{-\alpha}|\eta|^{-\beta} |\widehat{\psi}(R(\xi + \eta))| d\xi d\eta  $$
$$ \leq R^d \left(\iint |\widehat{f}(\xi)|^2 |\xi|^{-2\alpha} |\widehat{\psi}(R(\xi + \eta))| d\xi d\eta \right)^{1/2} \left(\iint\limits_{E} |\widehat{g}(\eta)|^2 |\eta|^{-2\beta} |\widehat{\psi}(R(\xi + \eta))| d\xi d\eta \right)^{1/2}$$
$$ \lesssim \left(\int |\widehat{f}(\xi)|^2 |\xi|^{-2\alpha} d\xi \right)^{1/2} \left(\int |\widehat{g}(\eta)|^2 |\eta|^{-2\beta} d\eta \right)^{1/2} $$
In the last step we used the following inequality
$$ R^d \int |\widehat{\psi}(R(\xi + \eta))| d\xi \leq \| \widehat{\psi} \|_1 \lesssim 1 $$
and the corresponding one where the integration is in $\eta$.

\vskip.125in 

\section{Proof of Theorem \ref{miracle2}} 
\label{miraclesection2} 

In this section we prove Theorem \ref{miracle2}. In order to study the integrability of
$$ \tilde{T}_{\mu}(f,g)(x)=\iint f(x+u-v) g(x+u) d\mu(u,v)$$
we look at
$$ \int \left|\tilde{T}_{\mu}(f,g)(x)\psi\left(\frac{x}{R}\right)\right| dx $$
where $\psi$ is a non-negative bump function with support in $[-2,2]$ and identically equal to $1$ on $[-1,1]$. We obtain bounds independent on $R$ and thus Theorem \ref{miracle2} follows by a standard limiting argument.

Assum $f$ and $g$ are positive functions so that we can write
$$ \int T_{\mu}(f,g)(x)\psi\left(\frac{x}{R}\right) dx .$$
Viewing the above expression from the Fourier side we obtain
$$ R^d \iint \widehat{f}(\xi)\widehat{g}(\eta) \widehat{\mu}(-\xi-\eta,\xi)\widehat{\psi}(R(\xi + \eta)) d\xi d\eta .$$
By using estimate \ref{decay2} we can bound the above expression by
$$ R^d \iint |\widehat{f}(\xi)| |\widehat{g}(\eta)|C(|\xi+\eta|+|\xi|)^{-\gamma}|\widehat{\psi}(R(\xi + \eta))| d\xi d\eta $$
$$ \lesssim R^d \iint |\widehat{f}(\xi)| |\widehat{g}(\eta)| |\xi|^{-\alpha}|\eta|^{-\beta} |\widehat{\psi}(R(\xi + \eta))| d\xi d\eta  $$
$$ \leq R^d \left(\iint |\widehat{f}(\xi)|^2 |\xi|^{-2\alpha} |\widehat{\psi}(R(\xi + \eta))| d\xi d\eta \right)^{1/2} \left(\iint\limits_{E} |\widehat{g}(\eta)|^2 |\eta|^{-2\beta} |\widehat{\psi}(R(\xi + \eta))| d\xi d\eta \right)^{1/2}$$
$$ \lesssim \left(\int |\widehat{f}(\xi)|^2 |\xi|^{-2\alpha} d\xi \right)^{1/2} \left(\int |\widehat{g}(\eta)|^2 |\eta|^{-2\beta} d\eta \right)^{1/2} $$
In the last step we used the following inequality
$$ R^d \int |\widehat{\psi}(R(\xi + \eta))| d\xi \leq \| \widehat{\psi} \|_1 \lesssim 1 $$
and the corresponding one where the integration is in $\eta$.

\vskip.125in 

\section{Proof of Theorem \ref{notsharp}} 

\vskip.125in 

Let $E_n$ denote the $n^{-\frac{1}{s}}$-neighborhood of
$$ P_n = \frac{1}{n^{\frac{1}{d}}}\left(\mathbb{Z}^d \cap [0,n^{1/d}]^d \right) $$
where $0<s<d/2$. It is known that the Hausdorff dimension of
$$ E = \bigcap\limits_{k=K}^{\infty}E_{2^{d \cdot 2^k}} , $$
where $K$ is a non-negative integer that we can choose, is $s$ and furthermore that it is Ahlfors-David regular. See, for example, \cite{Fal86}, \cite{Falc86}. See also \cite{M95} and \cite{W04} for a thorough description of the background material pertaining to fractal geometry and its connections with harmonic analysis.

Let $\mu_s$ be the $s$-dimensional Hausdorff measure on $E$ and take $\epsilon = n^{-\frac{1}{s}}$. In order to show that $\frac{\pi}{2}$ is not equitably represented in $\mathcal{A}(E)$ it is sufficient to establish
\begin{equation}\label{ninety}
\mu_s \times \mu_s \times \mu_s \left\{(x,y,z)\in\mathbb{R}^{3d}: - n^{-\frac{1}{s}} \leq \frac{x-z}{|x-z|} \cdot \frac{y-z}{|y-z|} \leq n^{-\frac{1}{s}} \right\} \gtrsim n^{-\frac{1}{s}} .
\end{equation}
The left hand side of \ref{ninety} can be bounded below by
\begin{equation}\label{ninetyE}
\mu_s \times \mu_s \times \mu_s \left\{(x,y,z)\in E_n^{3}: - n^{-\frac{1}{s}} \leq \frac{x-z}{|x-z|} \cdot \frac{y-z}{|y-z|} \leq n^{-\frac{1}{s}} \right\}
\end{equation}
for all $n$. Through some straight forward estimates we can see that if we have $(x_0, y_0, z_0) \in P_n^3$ such that $\theta(x_0, y_0, z_0) = \frac{\pi}{2}$ then
$$ B(x_0, n^{-\frac{1}{s}}) \times B(y_0, n^{-\frac{1}{s}})  \times B(z_0, n^{-\frac{1}{s}}) \subseteq \left\{(x,y,z)\in E_n^{3}: - n^{-\frac{1}{s}} \leq \frac{x-z}{|x-z|} \cdot \frac{y-z}{|y-z|} \leq n^{-\frac{1}{s}} \right\}.  $$
For the estimates we use the fact that $n^{-\frac{1}{s}} \ll n^{-\frac{1}{d}}$, given $n$ large enough, since $0<s<\frac{d}{2}$, and that the smallest distance between any two distinct points in $P_n^3$ is $n^{-\frac{1}{d}}$. This last observation also tells us that the sets
$$ B(x_0, n^{-\frac{1}{s}}) \times B(y_0, n^{-\frac{1}{s}})  \times B(z_0, n^{-\frac{1}{s}}) ,$$
where $B(x,r)$ denotes the ball in $\mathbb{R}^d$ with center $x$ and radius $r$, are disjoint for any two different choices of $(x_0, y_0, z_0) \in P_n^3$. Those sets are also finitely many so we see that we can bound \ref{ninetyE} below by
\begin{equation}\label{almostThere}
\left(\inf\limits_{x_0\in P_n}\mu_s(B(x_0,n^{-\frac{1}{s}}))\right)^3 \cdot \text{ the number of right angles in }{\mathcal A}(P_n).
\end{equation}
Note that a priori that $x_0$ need not be in in $E$. However since these inequalities hold for all $n$ we can choose to use $n = 2^{d \cdot 2^K}$, where we have chosen $K$ to be large enough. It is clear by our construction that $P_{2^{d \cdot 2^K}} \subseteq P_{2^{d \cdot 2^k}} $ for all $k \geq K$ and thus by our construction of $E$ we have $P_{2^{d \cdot 2^K}} \subseteq E$. Thus we can guarantee that the $x_0$ above is in $E$.

Recall that since $E$ is Ahlfors-David regular we know there exists a constant $C$ such that for all $x\in E$ and all $0<r\leq 1$ we have
$$ C^{-1}r^s \leq \mu_s(B(x,r)) \leq Cr^s .$$
Since we can make sure that $x_0$ above is in $E$ then this in particular means that we can bound \ref{almostThere} below by
$$ n^{-3} \cdot \text{ the number of right angles in }{\mathcal A}(P_n). $$
Using theorem \ref{counting} (scaling does not change number of right angles) the above is bounded below by
$$ n^{-3} \cdot n^{3-\frac{2}{d}} $$
and since $ s < \frac{d}{2} $ we have
$$ n^{-3} \cdot n^{3-\frac{2}{d}} > n^{-\frac{1}{s}} = \epsilon . $$
This shows that \ref{ninety} holds true.

\vskip.125in 

\section{Proof of Theorem \ref{counting}} 

\vskip.125in 

This is a relatively straight forward generalization of the argument of Apfelbaum and Sharir in \cite{AS05}. Assume for simplicity that $n$ is a $d$-th power and a multiple of $5$ so that all the quantities in the proof are integers. For a fixed $d$ then this assumption does not change the order of magnitude of the lower bound.

Recall that we write $f = O(g)$ if there exists an $x_0$ such that $f(x) \lesssim g(x)$ for all $x>x_0$. We write $f = \Omega(g)$ if there exists an $x_0$ such that $f(x) \gtrsim g(x)$ for all $x>x_0$. Finally we write $f=\Theta(g)$ if $f=O(g)$ and $f=\Omega(g)$.

Let $Q = \lbrace \frac{2}{5}n^{1/d} + 1, \ldots , \frac{3}{5}n^{1/d}\rbrace^d$ be the middle $\frac{1}{5}n^{1/d}\times \frac{1}{5}n^{1/d} \times \ldots \times \frac{1}{5}n^{1/d}$ portion of $P$. We have $|Q| = \frac{n}{5^d} = \Theta(n)$. For each pair of points in $Q$, the square of the distance between them is an integer of magnitude at most $\frac{d}{5^d}n^{2/d}$. Hence there are at most $\frac{d}{5^d}n^{2/d} = O(n^{2/d})$ distinct distances between the points of $Q$. For every point $x \in Q$ we take the spheres centered at $x$ and containing at least one point $p\in Q$. There are $O(n^{2/d})$ such spheres. Do this for all points in $Q$ and let $S$ denote the resulting set of such spheres. We thus have $|S| = O(n^{1+\frac{2}{d}})$. By choosing $Q$ to be small enough in $P$ we are guaranteed that for every point $p \in P$ on a sphere $\sigma \in S$, the point on $\sigma$, antipodal to $p$, is also in $P$.

For each $\sigma \in S$, let $m_{\sigma} = |P\cap\sigma |$ denote the number of lattice points on $\sigma$. We observe that $\sum_{\sigma\in S}m_{\sigma} \geq 2\binom{|Q|}{2}= \Omega(n^2)$, since in the sum we count every pair $p, p' \in Q$ exactly twice - once with $p$ at the center of the sphere and $p'$ on the sphere itself, and once the other way around. In a similar manner $\sum_{\sigma \in S}m_{\sigma} \leq |Q|\cdot |P| = O(n^2)$, so this sum is $\Theta(n^2)$. Let $\sigma \in S$ be one of the spheres and let $p,q,r\in\sigma\cap P$ be three distinct points such that $p$ and $r$ are antipodal points of $\sigma$. Then $\angle pqr = \pi/2$. There are $m_{\sigma}/2$ choices of an antipodal pair $p,r \in \sigma \cap P$ and $m_{\sigma}-2$ choices of a third point $q$. This yields $m_{\sigma}(m_{\sigma}-2)/2$ right angles on $\sigma$. The lower bound on the number of right angles in $P$ is obtained by summing over all the spheres of $S$. Note that each pair of points can be antipodal on at most one sphere, hence every angle is counted only once. This gives a lower bound of
$$ \frac{1}{2}\sum\limits_{\sigma \in S}m_{\sigma}(m_{\sigma}-2) \geq \frac{1}{2|S|}\left( \sum\limits_{\sigma\in S}m_{\sigma} \right)^2 - \sum\limits_{\sigma\in S}m_{\sigma} = \frac{1}{2|S|}\Theta(n^4) - \Theta(n^2), $$
where we have used the Cauchy-Schwarz inequality. Substituting $|S| = O(n^{1+\frac{2}{d}})$ in the inequality gives $\Omega(n^{3-\frac{2}{d}})$ right angles determined by the points of $P$.

\vskip.125in 

\section{Proof of Theorem \ref{notsharp2}} 

\vskip.125in 

\subsection{Construction of sets $E_n$} Define $n$ be the smallest integer greater than $R^{d-2}$, where $R$ is a square root of a large square free integer, and let $E_n$ denote the $n^{-\frac{1}{s}}$-neighborhood of 
$$ R^{-1} \{k \in {\mathbb Z}^d: |k|=R \}.$$ 

If $d \ge 5$, 
\begin{equation} \label{landau} \# \{k \in {\mathbb Z}^d: |k|=R \} \approx R^{d-2}. \end{equation}

In dimension four this is not necessarily the case, but does hold if $R^2$ is not divisible by $4$. See \cite{Mag07} and the references contained therein. 

We now establish (\ref{energy}). Note that in the language of Subsection \ref{discreteapplications}, we need to prove that the set $R^{-1} \{k \in {\mathbb Z}^d: |k|=R \}$ is $s$-adaptable for $0<s<\frac{d-1}{2}$. Let $\psi$ be a smooth cut-off supported in the unit ball. Define 
$$ d\mu_n(x)=n^{-1} n^{\frac{d-1}{s}} \sum_{\{k \in {\mathbb Z}^d: |k|=R \}} \phi(n^{\frac{1}{s}}(x-k/R))d\sigma(x),$$ where $d\sigma$ is the Leray measure on the unit sphere. 

It follows that 
$$ \int \int {|x-y|}^{-s} d\mu_n(x) d\mu_n(y)$$
$$={(n^{-1} n^{\frac{d-1}{s}})}^2 \sum_{\{k,l \in {\mathbb Z}^d: |k|=|l|=R \}} \int \int {|x-y|}^{-s} \phi(n^{\frac{1}{s}}(x-k/R))\phi(n^{\frac{1}{s}}(y-l/R)) d\sigma(x) d\sigma(y)$$
$$=I+II,$$ where 
$$ I={(n^{-1} n^{\frac{d-1}{s}})}^2 \sum_{\{k \in {\mathbb Z}^d: |k|=R \}} \int \int {|x-y|}^{-s} 
\phi(n^{\frac{1}{s}}(x-k/R)) \phi(n^{\frac{1}{s}}(y-k/R)) d\sigma(x)d\sigma(y)$$ and 
$$ II \approx {(n^{-1} n^{\frac{d-1}{s}})}^2 \sum_{\{k\not=l \in {\mathbb Z}^d: |k|=|l|=R \}} R^s {|k-l|}^{-s} \int \int \phi(n^{\frac{1}{s}}(x-k/R)) \phi(n^{\frac{1}{s}}(y-k/R)) d\sigma(x)d\sigma(y).$$

We have 
$$ II \approx n^{-2} R^s \sum_{\{k \not=l \in {\mathbb Z}^d: |k|=|l|=R \}} {|k-l|}^{-s}$$
\begin{equation} \label{erection2} \lesssim n^{-2} R^s \sum_{2^j \leq R} \sum_{|k|=R} 2^{-js} w_j(k),\end{equation} where 
$$ w_j(k)=\sum_{\{l \in {\mathbb Z}^d: |l|=R; 2^j \leq |k-l| \leq 2^{j+1} \}} 1.$$

We need the following estimate that follows from the main result in \cite{Mag07}. 
\begin{lemma} With the notation above, 
$$ |w_j(k)| \lesssim 2^{j(d-2)}.$$ 
\end{lemma} 

In view of the lemma, the expression in (\ref{erection2}) is 
$$ \lesssim n^{-2} R^s \sum_{2^j \leq R} \sum_{|k|=R} 2^{j(d-2-s)}$$
$$ \lesssim n^{-1} R^s \sum_{2^j \leq R} 2^{j(d-2-s)} \lesssim n^{-1} R^s R^{d-2-s}$$
$$ \lesssim n^{-1} R^{d-2} \lesssim 1.$$ 

This proves that $II \lesssim 1$. Since $I \lesssim 1$ by a direct calculation, the proof of (\ref{energy}) is complete. 

\subsection{Estimation of the size of ${\mathcal A}(E_n)$} It is not difficult to see that 
\begin{equation} \label{elephanterection} |{\mathcal A}(E_n)| \approx n^{-\frac{1}{s}} \cdot \# \{\theta(k^1,k^2,k^3): k^j \in \{m \in {\Bbb Z}^d: |m|=R \}.\end{equation}

We claim that 
$$ \# \{\theta(k^1,k^2,k^3): k^j \in \{m \in {\Bbb Z}^d: |m|=R \} \lesssim R^2.$$ 

To see this, observe that it is sufficient to count dot products of the form 
$$ (k^1-k^2) \cdot (k^1-k^3).$$ 

These are integers contained in $[-R^2, R^2]$ and there can be at most $2R^2$ such numbers. Plugging this into (\ref{elephanterection}), we see that 
$$ |{\mathcal A}(E_n)| \lesssim n^{-\frac{1}{s}} \cdot n^{\frac{2}{d-2}}$$ and the right hand side goes to $0$ as $n \to \infty$ if $s<\frac{d-2}{2}$. 

\vskip.125in 

\section{Proof of Theorem \ref{discrete}} 

\vskip.125in 

Since $P$ is $s$-adaptable we can thicken it into a set $E$, which is uniformly $s$-dimensional. Let $(x_0, y_0, z_0) \in P^3$ such that $\theta(x_0, y_0, z_0) = \theta_0$. Then straightforward estimates show that
$$ B\left(x_0, \frac{1}{8}N^{-\frac{1}{s}}\right) \times B\left(y_0, \frac{1}{8}N^{-\frac{1}{s}}\right) \times B\left(z_0, \frac{1}{8}N^{-\frac{1}{s}}\right)  ,$$
where $B(x,r)$ is the ball in $\mathbb{R}^d$ centered at $x$ with radius $r$, is contained in
$$ \{(x,y,z)\in E \times E \times E: \theta_0-N^{-\frac{1}{s}} \leq \theta(x,y,z) \leq \theta_0+N^{-\frac{1}{s}} \} .$$
Furthermore, since $P$ is $N^{-\frac{1}{s}}$ separated then two such sets, for two different $(x_0, y_0, z_0) \in P^3$, are disjoint. Finally note that
$$ \mu_P^s \left( B(x_0, \frac{1}{8}N^{-\frac{1}{s}}) \right) \leq N^{-1}. $$
Thus we can bound
$$ \# \{(x,y,z) \in P \times P \times P: \theta(x,y,z)=\theta_0 \} $$ 
above by
$$N^3 \cdot \mu_P^s \times  \mu_P^s \times  \mu_P^s \left\{(x,y,z)\in E^3: \theta_0-N^{-\frac{1}{s}} \leq \theta(x,y,z) \leq \theta_0+N^{-\frac{1}{s}} \right\} .$$
However by Theorem \ref{main} we know that
$$ \mu_P^s \times  \mu_P^s \times  \mu_P^s \left\{(x,y,z)\in E^3: \theta_0-N^{-\frac{1}{s}} \leq \theta(x,y,z) \leq \theta_0+N^{-\frac{1}{s}} \right\} \lesssim N^{-\frac{1}{s}} .$$
Thus we have shown
$$ \# \{(x,y,z) \in P \times P \times P: \theta(x,y,z)=\theta_0 \} \lesssim N^{3-\frac{1}{s}} $$
as claimed.

\end{document}